\documentclass{amsproc}

\newtheorem{Theorem}{Theorem}

\newtheorem{theorem}{Theorem}[section]
\newtheorem{lemma}[theorem]{Lemma}
\newtheorem{proposition}[theorem]{Proposition}

\theoremstyle{definition}

\theoremstyle{remark}
\newtheorem{remark}[theorem]{Remark}

\numberwithin{equation}{section}

\newcommand{\abs}[1]{\lvert#1\rvert}
\newcommand{\Abs}[1]{\Vert\lvert#1\rvert\Vert}

\newcommand{\R}{{\mathbb R}}
\newcommand{\Q}{{\mathbb Q}}
\newcommand{\Z}{{\mathbb Z}}

\newcommand{\BB}{{\mathcal B}}

\newcommand{\N}{{\mathbb N}}

\newcommand{\Fix}{{\rm Fix}}

\begin{document}

\title[Orbit equivalence types of circle diffeomorphisms]
{Orbit equivalence types of circle diffeomorphisms 
with a Liouville rotation number}


\author{Shigenori Matsumoto}
\address{Department of Mathematics, College of
Science and Technology, Nihon University, 1-8-14 Kanda, Surugadai,
Chiyoda-ku, Tokyo, 101-8308 Japan
}
\email{matsumo@math.cst.nihon-u.ac.jp
}
\thanks{The author is partially supported by Grant-in-Aid for
Scientific Research (C) No.\ 20540096.}
\subjclass{Primary 37E10,
secondary 37E45.}

\keywords{circle diffeomorphism, rotation number, Liouville number,
orbit equivalence, fast
approximation
by conjugation}

\maketitle

\date{\today }
\begin{abstract}
This paper is concerned about the orbit equivalence 
types of $C^\infty$ diffeomorphisms
of $S^1$ seen as nonsingular automorphisms of $(S^1,m)$, where
$m$ is the Lebesgue measure.
Given any Liouville number $\alpha$, it is shown that each of
the subspace
formed by type ${\rm II}_1$, ${\rm II}_\infty$, ${\rm III}_\lambda$
($\lambda>1$), ${\rm III}_\infty$ and ${\rm III}_0$ diffeomorphisms
are $C^\infty$-dense in
the space of the orientation preserving $C^\infty$ diffeomorphisms 
with rotation number $\alpha$.
\end{abstract}

\section{Introduction}
Let $(X,\BB,\mu)$ and $(X',\BB',\mu')$ be Lebesgue measure spaces.
A map $$f:(X,\BB,\mu)\to(X',\BB',\mu')$$
is called a {\em nonsingular isomorphism} (a nonsingular {\em automorphism}
when $(X',\BB',\mu')=(X,\BB,\mu)$), if $f$ is a bimeasurable bijection 
from a conull set of $X$ onto a conull set of $X'$ and if $f_*\mu$
is equivalent to $\mu'$.

Two nonsingular automorphisms
$$f:(X,\BB,\mu)\to(X,\BB,\mu)\ \mbox{ and }\
f':(X',\BB',\mu')\to(X',\BB',\mu')$$
are called {\em orbit equivalent} if there is
a measurable isomorphism
$$h:(X,\BB,\mu)\to(X',\BB',\mu')$$
which sends the $f$-orbit of almost any point to
an $f'$-orbit (regardless of the orders of the orbits).

A nonsingular automorphism
$$f:(X,\BB,\mu)\to(X,\BB,\mu)$$
is called {\em ergodic} if any $f$-invariant measurable subset of $X$
is either null or conull.

An ergodic nonsingular automorphism $f$ is said to be of type II if
there is a $f$-invariant
$\sigma$-finite measure $\nu$ which is equivalent to $\mu$. More
specifically $f$ is called of type ${\rm II}_1$ (resp.\ of type
${\rm II}_\infty$) if the $f$-invariant
measure $\nu$ is  finite (resp.\ infinite).
It is known that they are mutually exclusive.

An ergodic nonsingular automorphism $f$ is said to be of type III if
it is not of type  II. Type III nonsingular automorphism $f$ is
further classified according to the {\em ratio set} $r(f) 
\subset [0,\infty)$, defined as follows.
A number $\alpha\in[0,\infty)$ belongs to $r(f)$ if and only if
for any $A\in\BB$ with $\mu(A)>0$ and $\epsilon>0$, there is
$B\in\BB$ and $i\in\Z$ such that  $\mu(B)>0$, $B\subset A$, 
$f^i(B)\subset A$, and for any $x\in B$, the Radon-Nikodym derivative
satisies:
$$
\frac{d(f^{-i})_*\mu}{d\mu}(x)\in(\alpha-\epsilon,\alpha+\epsilon).
$$

The definition of the ratio set $r(f)$ does not depend on the choice of a measure $\mu$
from its equivalence class $[\mu]$.
Moreover it is an invariant of an orbit equivalence class.
The ratio set $r(f)$ is a closed subset of $[0,\infty)$, and
$r(f)\cap(0,\infty)$ is a multiplicative subgroup of $(0,\infty)$.
It is well known, easy to show, that $f$ is of type II if and only
if $r(f)=\{1\}$. For $f$ of type III, we have the following
classification.

\begin{itemize}
\item
$f$ is of type ${\rm III}_\lambda$ for some $\lambda>1$ if
$r(f)=\lambda^\Z\cup\{0\}.$

\item
$f$ is of type ${\rm III}_\infty$ if $r(f)=[0,\infty)$.

\item
$f$ is of type ${\rm III}_0$ if $r(f) = \{0,1\}$.
\end{itemize}

It is known (\cite{Kr}) that the set of ergodic nonsingular transformations
of type ${\rm II}_1$, ${\rm II}_\infty$, ${\rm III}_\lambda$,
${\rm III}_\infty$ each consists of one orbit equivalence class. But
the set of the ergodic nonsingular transformations of type ${\rm III}_0$
consists of various classes.

Let $F$ be the group
of the orientation preserving $C^\infty$ diffeomorphisms of $S^1$.
For $\alpha\in\R/\Z$, denote by $F_\alpha$ the subset of $F$
consisiting of those diffeomorphisms $f$ whose rotation number
$\rho(f)$ is $\alpha$.
If $\alpha$ is irrational, then any element $f\in F_\alpha$ is 
ergodic with respect to the Lebesgue measure $m$
(\cite{H}, p.86, \cite{KH} 12.7).

In \cite{Ka}, Y. Katznelson has shown that the diffeomorphisms
of any type raised above are $C^\infty$ dense in the union
of $F_\alpha$ for irrational $\alpha$.

In this paper we focus on a single subset $F_\alpha$. For any
$f\in F_\alpha$, $\alpha$ irrational, there is a unique
homeomorphism, denoted by $H_f$, such that $H_f(0)=0$
and $f=H_f R_\alpha H_f^{-1}$, where $R_\alpha$
stands for the rotation by $\alpha$.
Thus the unique $f$ invariant measure on $S^1$ is given by
$(H_f)_*m$. This implies
that either $(H_f)_*m$ is equivalent to $m$ or else singular
to $m$.

For a non Liouville number $\alpha$, it is shown (\cite{Y1}) that
$H_f$ is a $C^\infty$ diffeomorphism for any $f\in F_\alpha$.
That is,  $(H_f)_*m$ is equivalent to
the Lebesgue measure $m$, and hence any $f\in F_\alpha$ is of
type ${\rm II}_1$. 

But for a Liouville number $\alpha$, things
are quite different. There are $f\in F_\alpha$ for which the
unique $f$ invariant measure $(H_f)_*m$ is singular to $m$
(\cite{M1,M2}).
Such $f$ can never be of type ${\rm II}_1$.
Denote by $F_\alpha(\rm{T})$ the subset of $F_\alpha$ consisting of
diffeomorphisms of type T. The main result of this
paper is the following.

\begin{Theorem} \label{main}
For any Liouville number $\alpha$, each of the subsets $F_\alpha({\rm II}_1)$,
$F_\alpha({\rm II}_\infty)$, $F_\alpha({\rm III}_\lambda)$ for any $\lambda>1$,
$F_\alpha({\rm III}_\infty)$ and $F_\alpha({\rm III}_0)$
forms a $C^\infty$-dense subset of $F_\alpha$.
\end{Theorem}

The key fact for the proof is the result in
\cite{Y2} which states that even for a Liouville number $\alpha$,
the subset of elements $f\in F_\alpha$ such that
$H_f$ are $C^\infty$ diffeomorphisms is $C^\infty$ dense 
in $F_\alpha$.
Since the $C^\infty$ closure of our subset $F_\alpha({\rm T})$ is invariant
by the conjugation by an element $f\in F$, it suffices to show
the following proposition.

\begin{proposition} \label{p1}
For any $r\in\N$, there is an element $f\in F_\alpha({\rm T})$
such that $d_r(f,R_\alpha)<2^{-r}$, where $d_r$ is the $C^r$ 
distance and $\rm T=II_\infty,\ III_\lambda,\ III_\infty,\ III_0$.
\end{proposition}

Proposition \ref{p1} is proved by the method of
fast approximation by conjugacy with estimate, developed in \cite{FS}.
This is a qualitatively refined version of the method of successive
approximations originated by D. Anosov and A. Katok \cite{AK}
in the early 70's.

\bigskip\noindent
{\sc Acknowlegement}: The author expresses his gratitude to Masaki
Izumi, who brought his attention to this problem.

\section{Fast approximation method}

We assume throughout that $0<\alpha<1$ is
a Liouville number, i.\ e.\ for any $N\in\N$ and $\epsilon>0$,
there is $p/q$ ($p,q\in\N$, $(p,q)=1$)
such that
\begin{equation}\label{e11}
\abs{\alpha-p/q}<\epsilon q^{-N}.
\end{equation}

To prove Proposition \ref{p1} for $F_\alpha({\rm T})$, 
we will actually show the next
proposition.

\begin{proposition}\label{p2}
For any $r\in\N$, there are sequences $\alpha_n=p_n/q_n\in\Q\cap(0,1)$,
$(p_n,q_n)=1$,
 and $h_n\in F$ ($n\in\N$) such that the following {\rm (i)
 $\sim$
 (v)} hold.
Define $H_0={\rm Id}$, $f_0=R_{\alpha_1}$, and for any $n\in\N$
$$
H_n=h_1\cdots h_n\ \mbox{ and }\ f_n=H_n R_{\alpha_{n+1}}H_n^{-1}.$$

\noindent
{\rm (i)} $\alpha_n\to\alpha$.
\\
{\rm (ii)} $R_{\alpha_n}$ commutes with $h_n$. 
\\
{\rm (iii)} $H_n$ converges uniformly to a homemorphism $H$. 
\\
{\rm (iv)}
$$\abs{\alpha-\alpha_1}<2^{-r-1},\ \mbox{and}\ \
d_{n+r}(f_{n-1},f_{n})<2^{-n-r-1},\ \ \forall n\geq1.$$
{\rm (v)} The limit $f$ of $\{f_n\}$ is an element of $F_\alpha({\rm T})$.
\end{proposition}

Notice that (iv) implies that the limit $f$ of $f_n$
is a $C^\infty$ diffeomorphism
such that
$d_r(f,R_\alpha)<2^{-r}$.

Condition (ii) is useful to establish (iv), since
then
\begin{eqnarray*}
f_{n-1}-f_n&=&H_nR_{\alpha_n}H_n^{-1}-H_nR_{\alpha_{n+1}}H_n^{-1},\\
f_{n-1}^{-1}-f_n^{-1}&=&H_nR_{-\alpha_n}H_n^{-1}-H_nR_{-\alpha_{n+1}}H_n^{-1},
\end{eqnarray*}
and these can be estimated using Lemma \ref{l2} below.

Next we shall summerize inequalities needed to establish Proposition
\ref{p2}. All we need are polynomial type estimates
whose degree and coefficients can be arbitrarily large.
The inequalities below are sometimes far from being optimal.

For a $C^\infty$ function $\varphi$ on $S^1$, we define 
as usual the $C^r$ norm
$\Vert\varphi\Vert_r$ ($0\leq r<\infty$) by
$$
\Vert\varphi\Vert_r=\max_{0\leq i\leq r}\sup_{x\in S^1}\abs{\varphi^{(i)}(x)}.$$

For $f,g\in F$, define
\begin{eqnarray*}
\Abs{f}_r&=&\max\{\Vert f-{\rm id}\Vert_r,\ \Vert f^{-1}-{\rm id}\Vert_r,1\},\\
d_r(f,g)&=&\max\{\Vert f-g\Vert_r,\ \Vert f^{-1}-g^{-1}\Vert_r\}.
\end{eqnarray*}

The term $\Abs{f}_r$ is used to show that $f$ is not so large in the
$C^r$-topology. 
We have included $1$ in its definition 
because then it becomes possible to reduce inequalities
from the Fa\`a di Bruno formula (\cite{H}, p.42 or \cite{S})
by virtue of the following;
$$ \Abs{f}_r^i\leq\Abs{f}_r^r\ \ {\rm if}\ \ i\leq r.$$
On the other hand $d_r(f,g)$ is useful for showing
$f$ and $g$ are near in the $C^r$-topology. 
We get the following inequality from the Fa\`a di Bruno formula.

Below we denote by C(r) an arbitrary constant which depends only
on $r$.

\begin{lemma} \label{l1}
For $f,g\in F$ we have
\begin{eqnarray*}
\Vert fg-g\Vert_r&\leq& C(r)\Vert f-{\rm Id}\Vert_r\,
\Abs{g}^r_r,\\
\Abs{fg}_r&\leq & C(r)\,\Abs{f}_r^r\,\Abs{g}_r^r.
\end{eqnarray*}

\qed
\end{lemma}

The next lemma can be found as Lemma 5.6 of \cite{FS} or
as Lemma 3.2 of \cite{S}.

\begin{lemma} \label{l2}
For $H\in F$ and $\alpha,\beta\in\R/\Z$,
$$
d_r(HR_\alpha H^{-1},HR_\beta H^{-1})\leq
C(r)\,\Abs{H}_{r+1}^{r+1}\,\abs{\alpha-\beta}.
$$
\qed
\end{lemma}

For $Q\in\N$, denote by $\pi_Q:S^1\to S^1$ the cyclic $Q$-fold covering map.

\begin{lemma} \label{l3}
Let $h$ be a lift of $\hat h\in F$ by $\pi_Q$ and assume $\Fix(h)\neq\emptyset$. 
Then we have for any $r\geq0$
\begin{eqnarray*}
\Vert h-{\rm Id}\Vert_r&=&\Vert \hat h-{\rm Id}\Vert_r\,Q^{r-1}, \\
\Abs{h}_r&\leq&\Abs{\hat h}_r\,Q^{r-1}.
\end{eqnarray*}
\end{lemma}

{\sc Proof}. Just notice that a lift $\tilde h$ of $h$ to $\R$ is 
the conjugate of a lift $\tilde {\hat h}$ of $\hat h$ by a homothety by $Q$,
i.\ e.\ $\tilde h(x)=Q^{-1}\tilde{\hat h}(Qx)$.
\qed

\bigskip
First of all we give
beforehand a sequence of diffeomorphisms $\hat h_n\in F$ such that
$\hat h_n(0)=0$.
The sequence $\{\hat h_n\}$ depends upon the type 
$\rm T=II_\infty,\ III_\lambda,\ III_\infty,\ III_0$,
and will be constructed concretely in the later sections.

Next we choose a
sequence of rationals $\alpha_n=p_n/q_n$ inductively
in a way to be explained
shortly,
 and
set $h_n$ to be the lift of $\hat h_n$ by the cyclic $Q_n$-fold covering
map $$\pi_{Q_n}:S^1\to S^1$$
such that ${\rm Fix}(h_n)\neq\emptyset$, where $Q_n=K(n)q_n$
and the integer $K(n)$ are chosen to depend only on 
$\{\hat h_i\}_{i\in\N}$ and
$q_1,\cdots,q_{n-1}$.
Notice that then  condition (ii) of Proposition \ref{p2} is 
automatically satisfied.

We always choose the rationals $\alpha_n$ so as to satisfy
$$
\abs{\alpha_{n+1}-\alpha}<\abs{\alpha_{n}-\alpha},\ \ \forall n\in\N.$$
Therfore we have
$$
\abs{\alpha_n-\alpha_{n+1}}\leq 2\abs{\alpha-\alpha_n}.$$

In the sequal, we denote any constant which
depends only on $r$, $\{\hat h_i\}_{i\in\N}$  and 
$\alpha_1,\cdots, \alpha_{n-1}$ by $C(n,r)$. 
Thus $C(n,r)$ depends only on the innitial data about
$\hat h_i$
and the previous step of the induction.
We also denote any positive
integer
which depends on $r$, $\{\hat h_i\}_{i\in\N}$ and 
$\alpha_1,\cdots, \alpha_{n-1}$ by $N(n,r)$.

By Lemma \ref{l3}, we have for any $1\leq i< n$,
\begin{equation}\label{e1}
\Abs{h_i}_{n+r+1}\leq \Abs
{\hat h_i}_{n+r+1}Q_i^{n+r}=C(n,r),
\end{equation}

and
\begin{equation}\label{e2}
\Abs{h_n}_{n+r+1}\leq \Abs
{\hat h_n}_{n+r+1}K(n)^{n+r}q_n^{n+r}=C(n,r)q_n^{N(n,r)}.
\end{equation}
Of course the two $C(n,r)$'s in (\ref{e1}) and (\ref{e2})
are different.
Now we obtain inductively using Lemma \ref{l1} that
\begin{equation}\label{e3}
\Abs{H_n}_{n+r+1}\leq C(n,r)q_n^{N(n,r)}.
\end{equation}
The terms $C(n,r)$ and $N(n,r)$ in (\ref{e3}) are computed
from (\ref{e1}) and (\ref{e2})
by applying Lemma \ref{l1} successively.

By Lemma \ref{l2} and (\ref{e3}),
\begin{eqnarray}\label{e4}
d_{n+r}(f_{n-1},f_n)
&=&d_{n+r}(H_nR_{\alpha_n}H_n^{-1},H_nR_{\alpha_{n+1}}H_n^{-1})
\\
&\leq&C(n,r)\,q_n^{N(n,r)}\abs{\alpha_n-\alpha_{n+1}}
\\
&\leq&C(n,r)\,q_n^{N(n,r)}\abs{\alpha-\alpha_n},\nonumber
\end{eqnarray}
for some other $C(n,r)$ and $N(n,r)$.

In order to obtain (iv) of Proposition \ref{p2}, 
the rational $\alpha_n=p/q$ have only to satisfy
$$C(n,r)q^{N(n,r)}\abs{\alpha-p/q}<2^{-n-r-1},
$$
that is,
\begin{equation}\label{e5}
\abs{\alpha-p/q}<2^{-n-r-1}C(n,r)^{-1}q^{-N(n,r)}.
\end{equation}

The terms
$$
\epsilon=2^{-n-r-1}C(n,r)^{-1}\ \mbox{ and }\
N=N(n,r)$$
are already determined beforehand or by the previous step of the induction.
Since $\alpha$ is Liouville, there exists a rational $p/q$
which satisfies (\ref{e11}) for these values of $\epsilon$
and $N$.
Setting it $p_n/q_n$, we establish (iv) for the
$n$-th step of the induction.

In fact there are infinitely many choices of $p_n/q_n$, 
which enables us to require more. First of all,
since
$$
\Vert H_n-H_{n-1}\Vert_0\leq{\rm Lip}(H_{n-1})\Vert h_n-{\rm Id}\Vert_0
\leq{\rm Lip}(H_{n-1})Q_n^{-1}\leq{\rm Lip}(H_{n-1})q_n^{-1},$$
we obtain that $H_n$ converges uniformly
to a continuous map $H$ simply if we choose $q_n$ large enough compared with 
the Lipschitz constant ${\rm Lip}(H_{n-1})$ of $H_{n-1}$.
It is easier to obtain that $H_n^{-1}$ converges. Then
$H$ is a homeomorphism, getting condition (iii).

There are other requirements for the choice of $\alpha_n$, which
will appear in the next section.

\section{Type ${\rm III}_\lambda$: the construction of $\hat h_n$}

In this and the next two sections, we shall prove Proposition \ref{p2}
for $T={\rm III}_\lambda$ and $\lambda>1$.
The diffeomorphism $\hat h_n$ ($n\in\N$) is constructed as follows.
First consider two affine maps of slope $\lambda^{1/2}$
and $\lambda^{-1/2}$, depicted in Figure 1. They intersect at
points $(0,0)$ and $(a,1-a)$.
Choose a rational number $\delta_n>0$ such that $\delta_n\downarrow 0$. 
(The precise condition will be given later.)
We join the two affine maps on the intervals $(-\delta_n,\delta_n)$
and $(a-\delta_n,a+\delta_n)$
 using bump functions. See Figure 2.
Finally the diffeomorphism $\hat h_n$ is obtained by adding a positive
number so that $\hat h_n$ has a fixed point at $0$.

\begin{figure}\caption{ }
\unitlength 0.1in
\begin{picture}( 39.4000, 33.6000)(  9.8000,-36.6000)
%
{%
\special{pn 8}%
\special{pa 1260 500}%
\special{pa 4890 500}%
\special{pa 4890 3640}%
\special{pa 1260 3640}%
\special{pa 1260 500}%
\special{pa 4890 500}%
\special{fp}%
}%
%
{%
\special{pn 8}%
\special{pa 1260 3630}%
\special{pa 4900 500}%
\special{fp}%
}%
%
{%
\special{pn 8}%
\special{pa 4310 3460}%
\special{pa 4920 300}%
\special{fp}%
}%
%
{%
\special{pn 8}%
\special{pa 4740 3070}%
\special{pa 980 3660}%
\special{fp}%
}%
\put(38.0000,-22.8000){\makebox(0,0)[lb]{slope $\lambda^{1/2}$}}%
\put(24.2000,-32.4000){\makebox(0,0)[lb]{slope $\lambda^{-1/2}$}}%
\end{picture}%
\end{figure}

\begin{figure}\caption{ }
\unitlength 0.1in
\begin{picture}( 47.5000, 23.0000)(  9.3000,-27.6000)
%
{{%
\special{pn 13}%
\special{pa 930 2500}%
\special{pa 2010 1000}%
\special{fp}%
}}%
%
{{%
\special{pn 13}%
\special{pa 2020 1010}%
\special{pa 3430 460}%
\special{fp}%
}}%
%
{{%
\special{pn 13}%
\special{pa 1760 1370}%
\special{pa 1782 1346}%
\special{pa 1802 1320}%
\special{pa 1824 1296}%
\special{pa 1844 1272}%
\special{pa 1910 1200}%
\special{pa 1932 1178}%
\special{pa 1956 1156}%
\special{pa 1978 1134}%
\special{pa 2004 1114}%
\special{pa 2028 1094}%
\special{pa 2080 1058}%
\special{pa 2108 1040}%
\special{pa 2134 1024}%
\special{pa 2162 1006}%
\special{pa 2214 970}%
\special{pa 2240 950}%
\special{pa 2264 930}%
\special{pa 2290 910}%
\special{fp}%
}}%
%
{{%
\special{pn 4}%
\special{pa 1730 1370}%
\special{pa 1730 2100}%
\special{fp}%
}}%
%
{{%
\special{pn 4}%
\special{pa 1990 1000}%
\special{pa 1990 2500}%
\special{fp}%
}}%
%
{{%
\special{pn 4}%
\special{pa 2260 910}%
\special{pa 2260 2070}%
\special{fp}%
}}%
%
{{%
\special{pn 4}%
\special{pa 2340 880}%
\special{pa 2340 2480}%
\special{fp}%
}}%
%
{{%
\special{pn 4}%
\special{pa 1670 1480}%
\special{pa 1670 2520}%
\special{fp}%
}}%
%
{{%
\special{pn 4}%
\special{pa 2620 780}%
\special{pa 2620 2490}%
\special{fp}%
}}%
\put(20.5000,-18.0000){\makebox(0,0)[lb]{$\delta_n$}}%
\put(23.8000,-20.5000){\makebox(0,0)[lb]{$\delta_n$}}%
\put(17.8000,-17.8000){\makebox(0,0)[lb]{$\delta_n$}}%
%
{{%
\special{pn 4}%
\special{pa 1410 1830}%
\special{pa 1410 2510}%
\special{fp}%
}}%
\put(14.6000,-20.9000){\makebox(0,0)[lb]{$\delta_n$}}%
\put(12.6000,-26.5000){\makebox(0,0)[lb]{$\partial^+\hat I_n^+$}}%
\put(24.8000,-26.5000){\makebox(0,0)[lb]{$\partial^-\hat I_n^-$}}%
\put(15.8000,-28.9000){\makebox(0,0)[lb]{$\partial^+\hat J_n^+$}}%
\put(22.4000,-28.7000){\makebox(0,0)[lb]{$\partial^-\hat J_n^-$}}%
%
\put(19.1000,-26.0000){\makebox(0,0)[lb]{}}%
\put(19.2000,-26.5000){\makebox(0,0)[lb]{$0$}}%
%
{{%
\special{pn 13}%
\special{pa 3420 2500}%
\special{pa 4640 2030}%
\special{fp}%
}}%
%
{{%
\special{pn 13}%
\special{pa 4630 2030}%
\special{pa 5680 510}%
\special{fp}%
}}%
%
{{%
\special{pn 13}%
\special{pa 4440 2100}%
\special{pa 4466 2080}%
\special{pa 4490 2060}%
\special{pa 4514 2038}%
\special{pa 4538 2018}%
\special{pa 4560 1996}%
\special{pa 4584 1974}%
\special{pa 4606 1950}%
\special{pa 4630 1928}%
\special{pa 4652 1906}%
\special{pa 4676 1884}%
\special{pa 4698 1862}%
\special{pa 4746 1818}%
\special{pa 4770 1798}%
\special{pa 4794 1776}%
\special{pa 4800 1770}%
\special{fp}%
}}%
%
{{%
\special{pn 4}%
\special{pa 4640 2040}%
\special{pa 4640 2470}%
\special{fp}%
}}%
%
{{%
\special{pn 4}%
\special{pa 4860 1720}%
\special{pa 4860 2430}%
\special{fp}%
}}%
%
{{%
\special{pn 4}%
\special{pa 4380 2140}%
\special{pa 4380 2460}%
\special{fp}%
}}%
%
\put(49.3000,-16.2000){\makebox(0,0)[lb]{}}%
%
{{%
\special{pn 4}%
\special{pa 4910 1620}%
\special{pa 4910 2580}%
\special{fp}%
}}%
%
{{%
\special{pn 4}%
\special{pa 5140 1320}%
\special{pa 5140 2590}%
\special{fp}%
}}%
%
{{%
\special{pn 4}%
\special{pa 4330 2150}%
\special{pa 4330 2600}%
\special{fp}%
}}%
%
{{%
\special{pn 4}%
\special{pa 4030 2260}%
\special{pa 4030 2590}%
\special{fp}%
}}%
\put(41.0000,-24.8000){\makebox(0,0)[lb]{$\delta_n$}}%
\put(44.2000,-22.7000){\makebox(0,0)[lb]{$\delta_n$}}%
\put(46.9000,-23.5000){\makebox(0,0)[lb]{$\delta_n$}}%
\put(49.7000,-22.0000){\makebox(0,0)[lb]{$\delta_n$}}%
\put(38.8000,-27.3000){\makebox(0,0)[lb]{$\partial^+\hat I_n^-$}}%
\put(42.4000,-28.9000){\makebox(0,0)[lb]{$\partial^+\hat J_n^-$}}%
\put(46.4000,-25.9000){\makebox(0,0)[lb]{$a$}}%
\put(48.1000,-28.6000){\makebox(0,0)[lb]{$\partial^-\hat I_n^-$}}%
\put(50.4000,-27.1000){\makebox(0,0)[lb]{$\partial^-\hat I_n^-$}}%
\end{picture}%
\end{figure}

Choose four {\em rational} points $\partial^+\hat J^+_n$, $\partial^-\hat J_n^-$,
$\partial^+\hat J_n^-$ and $\partial^-\hat J_n^+$ each near the points
$-\delta_n$, $\delta_n$, $a-\delta_n$ and $a+\delta_n$ as in Figure 
2.

Define two intervals by
$$
\hat J_n^-=[\partial^-\hat J_n^-,\partial^+\hat J_n^-],\ \ 
\hat J_n^+=[\partial^-\hat J_n^+,\partial^+\hat J_n^+].
$$
Then the diffeomorphism $\hat h_n$
is affine of slope $\lambda^{\pm1/2}$ on $\hat J_n^\pm$.
Next define subintervals $\hat I_n^\pm$ of $\hat J_n^\pm$ in
such a way that the connected components of $\hat J_n^\pm-\hat I_n^\pm$
have length $\delta_n$. Their boundary points $\partial^\pm I_n^\pm$
is denoted in Figure 2.

Let
\begin{equation}\label{e}
m(\hat h_n(\hat I_n^-\cup I_n^+))=1-\delta_n',
\end{equation}
where $m$ stands for the Lebesgue measure as before.
If one 
choose $\delta_n$ tending rapidly to $0$, and $\hat J_n^\pm$
appropriately, then we have
\begin{equation} \label{prod}
\prod_{n=1}^\infty (1-\delta_n')>9/10.
\end{equation}

Let $K'(n)$ be the least common multipliers of the denominators
of the eight rational points $\partial^\pm\hat J_n^\pm$ and
$\partial^\pm\hat I_n^\pm$. Define the number $K(n)$ 
inductively as follows. Let $K(1)=1$. When we are defining $K(n)$,
we already decided $q_{n-1}$. Set
$K(n)=q_{n-1}K'(n-1)$, and notice that this choice of $K(n)$
satisfies the condition of the previous section.

As before,
 set  $Q_n=K(n)q_n$ and $h_n$ 
to be the lift of $\hat h_n$ by the cyclic $Q_n$-fold covering
map $\pi_{Q_n}$
such that $0$ is a fixed point of $h_n$. Let
$$
J_n^{\pm}=\pi_{Q_n}^{-1}(\hat J_n^\pm),\ \ \
I_n^{\pm}=\pi_{Q_n}^{-1}(\hat I_n^\pm).$$
Now $J_n^{\pm}$ and $I_n^\pm$ each consists of $Q_n$ small intervals.
Each of these small intervals, as well as the fixed points of $h_n$,
are left fixed by $h_{n+1}$. Thus we have inductively
\begin{equation} \label{e31}
\mbox{The boundary points of $J_n^\pm$ and $I_n^\pm$ are fixed
by $h_m$ if $m>n$.}
\end{equation}

Finally we also assume the following. 
\begin{eqnarray} 
& q_{n+1}^{-1}<2^{-1} Q_n^{-1}\delta_n \label{e32}\\
& \abs{\alpha_n-\alpha}<\delta^2Q_n^{-2}q_n^{-1}\label{e33}.
\end{eqnarray}
One can assume (\ref{e32}) since $q_{n+1}$ can be chosen arbitrarily
large compared with the previous data, and (\ref{e33}) since this
takes the form (\ref{e11}) for the definition of Liouville numbers.

\section{Type ${\rm III}_\lambda$: the proof that $r(f)\subset \lambda^\Z\cup\{0\}$}

At this point we have already constructed the sequences $\{h_n\}$
and $\{\alpha_n\}$ in Proposition \ref{p2}. Thus for
$$
H_n=h_1\cdots h_n,\ \ f_n=H_n R_{\alpha_{n+1}}
H_n^{-1},
$$
we have obtained a $C^\infty$ diffeomorphism $f$ as the limit of
$\{f_n\}$ and a homeomorphism $H$ as the limit of $H_n$. 
They satisfy
$$
f=H R_\alpha H^{-1}.$$

What is left is to show that $f$ is a nonsingular automorphism of $(S^1,m)$
of type ${\rm III}_\lambda$. The purpose of this section is 
to show that $r(f)\subset \lambda^\Z\cup\{0\}$. For this,
it suffices to construct a Borel set $\Xi$ (in fact a Cantor set) 
of positive measure
which satisfies the following proposition.

\begin{proposition} \label{l41}
If $\xi\in\Xi$ and $f^i(\xi)\in\Xi$ for some $i\in\Z$, then
$(f^i)'(\xi)\in\lambda^\Z$.
\end{proposition}

Notice that the Radon-Nikodym derivative is just the usual derivative:
$$
\frac{d(f^{-i})_*m}{dm}(\xi)=(f^i)'(\xi).
$$

For $n\in \N$, let
$$
X_n=\bigcap_{j=1}^n(I_j^-\cup I_j^+),\ \ 
Y_n=\bigcap_{j=1}^n(J_j^-\cup J_j^+)
\ 
\mbox{ and }\ \
X=\bigcap_{j=1}^\infty(I_j^-\cup I_j^+).
$$

Setting
$$
\Xi=H(X),\ \ \ \Xi_n=H(X_n),$$
we shall show that $\Xi$ satisfies Proposition \ref{l41}.

Denote $H^{(n+1)}=\lim_{m\to\infty}h_{n+1} h_{n+2} \cdots
 h_m.$
Thus $H^{(n+1)}$ is a homeomorphism which satisfies
$$
H=H_n H^{(n+1)}.$$

For $x\in X_n$, denote by $[x]_n$ the connected component of
$x$ in $X_n$. For $x,\ x'\in X_n$, denote $x\sim_n x'$
if $[x]_n=[x']_n$.

\begin{lemma} \label{l43} Suppose $n<m$.

{\rm (1)} If $x,x'\in X_m$ satisfies $x\sim_m x'$,
then $x\sim_n y$.

{\rm (2)} For any $x\in X_n$, $h_{n+1}h_{n+2}\cdots h_{m}(x)\sim_n x$
 and
$H^{(n+1)}x\sim_n x$.
\end{lemma}

{\sc Proof}. (2) is a direct consequence of (\ref{e31}).
\qed

First of all, we have to show the following.

\begin{lemma} \label{l42}
The set $\Xi$ has positive Lebesgue measure.
\end{lemma}

{\sc Proof}.
Let us compute $m(\Xi_n)$. For $n=1$,
$\Xi_1=h_1(I^-_1\cup I_1^+)$, and clearly 
$$
m(\Xi)=m(\hat h_1(\hat I^-_1\cup \hat I_1^+))=1-\lambda'_1.$$
For $n=2$,
$$
\Xi_2=H((I_1^-\cup I_1^+)\cap(I_2^-\cup I_2^+)
=h_1((I_1^-\cup I_1^+)\cap h_2(I_2^-\cap I_2^+)),$$
by virtue of Lemma \ref{l43}.
Also by (\ref{e31}), the conditional probabilities of
$h_2(I_2^-\cap I_2^+)$ conditioned to $I_1^-$ and  $I_1^+$
coincide with $m(\hat h_2(\hat I_2^-\cup\hat I_2^+))$.
On the other hand $h_1$ is affine on $I_1^-$ and  $I_1^+$.
Therefore 
$$m(\Xi_2)=m(\hat h_1(\hat I_1^-\cup \hat I_1^+))
m(\hat h_2(\hat I_2^-\cup \hat I_2^+))=(1-\lambda'_1)(1-\lambda_2').
$$

Successive use of (\ref{e31}) enables us to conclude that in
general for $n\in\N$
$$
m(\Xi_n)=\prod_{j=1}^n (1-\lambda'_j),$$
and thus $m(\Xi)$ is positive by the assumption (\ref{prod}).
\qed

\begin{remark} \label{r41}
The above proof shows that any nonempty open subset
of $\Xi$ has positive measure. Furthermore any component of $\Xi_n\cap H(I_n^{+})$ has
the same measure, as well as any component of $\Xi_n\cap H(I_n^{-})$.
\end{remark}

\begin{lemma}\label{l44}
If $f^i(\xi)\in\Xi$ and $\abs{i}\leq n$, then
$f_n^i(\xi)\in H_n(Y_n)$.
\end{lemma}

{\sc Proof}. Let $x=H^{-1}(\xi)$ and $x_n=H_n^{-1}(\xi)=H^{(n+1)}(x)$.
By the assumption, the point $R_\alpha^i(x)=H^{-1}(f^i(\xi))$
belongs to $X\subset X_n$.
To show the lemma, it suffices to prove
$R_{\alpha_{n+1}}^i(x_n)=H_n^{-1}(f_n^{i}(\xi))$
is a point of $Y_n$. This follows once we show that
$
\abs{R_\alpha^i(x)-R^i_{\alpha_{n+1}}(x_n)}
$
is smaller than the width $Q_n^{-1}\delta_n$ of the
connected components of $Y_n\setminus X_n$.
Now we have
$$\abs{R^i_\alpha(x)-R^i_{\alpha_{n+1}}(x_n)}\leq
\abs{i}\abs{\alpha-\alpha_{n+1}}+\abs{x_n-x}.
$$
For $\abs{i}\leq n$,
$$\abs{i}\abs{\alpha-\alpha_{n+1}}\leq 
n\abs{\alpha-{\alpha_{n+1}}}<n\abs{\alpha-{\alpha_{n}}}
<n\delta_n^2 Q_n^{-2}q_n^{-1}<
2^{-1}Q_n^{-1}\delta_n
$$
by (\ref{e33}), and
$$\abs{x_n-x}<q_{n+1}^{-1}<2^{-1}Q_n^{-1}\delta_n
$$
by (\ref{e32}). This completes the proof.
\qed

\bigskip

\begin{lemma}\label{l45}
If $\xi,\ f^i(\xi)\in\Xi$ and $\abs{i}\leq n$, then
 $(f_n^i)'(\xi)\in\lambda^\Z$.
\end{lemma}

{\sc Proof}. Let $x=H^{-1}(\xi)$ and
$x_k=H_k^{-1}(\xi)$ for any $1\leq k\leq n$.
Then by Lemma \ref{l43}, since $x\in X\subset X_n$, we have
$x_k\in X_k$. Therefore $h_k'(x_k)=\lambda^{\pm 1/2}$.
Let $y_n=R_{\alpha_{n+1}}(x_n)$ and 
$y_k=h_{k+1}\cdots  h_n(y_n)$. From Lemma \ref{l44}, it follows that
 $y_n\in Y_n$, and
therefore by Lemma \ref{l43} applied to $\{Y_k\}$, we have $y_k\in Y_k$,
and in particular $h_k'(y_k)=\lambda^{\pm 1/2}$.
Now we have
$$
(f_n^i)'(\xi)=(h_1)'(y_1)\cdots (h_n)'(y_n)\,\cdot\,
(h_1)'(x_1)^{-1}\cdots (h_n)'(x_n)^{-1}\in\lambda^\Z.$$
\qed

The proof that $\Xi$ satisfies Proposition \ref{l41}
is complete by the following lemma.

\begin{lemma}\label{l46}
If $\xi,\ f^i(\xi)\in\Xi$, then for any large $n$,
 $(f_n)'(\xi)=f'(\xi)$.
\end{lemma}

{\sc Proof}. Here we shall give a proof which can also be applicable
when we show the denseness of type ${\rm III}_\infty$ diffeomorphisms
in later section.
Since $f^i_n$ converges to $f$
in $C^\infty$ topology and since $(f_n)'(\xi)-(f_{n+1})'(\xi)$
is either $1$ or $\lambda^{\pm1}$,  it must be 1 for any large $n$.
\qed

\section{Type ${\rm III}_\lambda$: the proof that $\lambda\in r(f)$}

Let $f$ and $\Xi$ be as before. Let $$f_\Xi:\Xi\to \Xi$$
be the first return map of $f$. As is well known, the ratio set
$r(f_\Xi)$ is the same as $r(f)$. 
To show that $\lambda\in r(f)$, we need to establish the
following proposition.

\begin{proposition} \label{p51}
For any closed subset $K\subset\Xi$ of positive measure, there
exist a point $\xi\in K$ and a number $i\in\Z$ such that 
$f^i(\xi)\in K$ and $(f^i)'(\xi)=\lambda$.
\end{proposition}

In fact this is enough for showing that $\lambda\in r(f_\Xi)$.
For, given any Borel subset $A\subset\Xi$, there is a closed
subset $K$ of positive measure
contained in the set of the points of density of $A$.
If there are $\xi$ and $i$ as in the proposition for this closed subset
$K$, then there is a small
interval $J_1$ centered at $\xi$ of radius $r$
such that $(f^i)'=\lambda$ on $J_1\cap \Xi\cap f^{-i}(\Xi)$.
Consider an interval $J_2$ centered at $f^i(\xi)$ of radius $(f^i)'(\xi)r$.
If $r$ is small enough, then we have
$$
m(J_1\cap A)>(1-\epsilon)m(J_1)\ \mbox{ and }\
m(J_2\cap A)>(1-\epsilon)m(J_2),$$
for some small
$\epsilon>0$, since $\xi$ and $f^i(\xi)$ is a point of density of $A$ .
Also if $r$ is small, $f^{-1}(J_2)$
almost coincides with $J_1$, and $f'$ is almost constant on $J_1$. 
 This shows that $B=(J_1\cap A)\cap f^{-1}(J_2\cap A)$ has
positive measure. For any $\eta\in B$, we have $f^i(\eta)\in A$
and $(f^i)'(\eta)=\lambda$, showing that $\lambda\in r(f_\Xi)$.

\bigskip

The rest of this section is devoted to the proof of Proposition
\ref{p51}. It is easier to pass from $\Xi$ to $X=H^{-1}(\Xi)$.
So let $\mu=H^{-1}_*m$, and choose once and for all an arbitrary
closed subset $C$ of $X$ of positive $\mu$ measure.
We shall show Proposition \ref{p51} for $K=H(C)$.

Our overall strategy is as follows. After choosing a large number $n$,
we shall construct points $x_k$, $y_k$ inductively for $k\geq n$.
They satisfy the following conditions, where $k$ is any number
bigger than $n$.

\begin{enumerate}
\item
$x_n$ (resp.\ $y_n$) is the midpoint of a component of
$X_n\cap I_n^-$ (resp.\ $X_n\cap I_n^+$), and $x_n\sim_{n-1}y_n$.
\item
$x_k$ and $y_k$ are the midpoints of  components of $X_k$ 
\item
 There is $i\in\Z$, independent of $k$, such that $y_k=R_{\alpha_k}^i(x_k)$.
\item
 $x_k\sim_{k-1}x_{k-1}$ and $y_k\sim_{k-1}y_{k-1}$.
\item
 $\mu([x_k]_k\cap C)>0$ and $\mu([y_k]_k\cap C)>0$.
\end{enumerate}

Let us show that this suffices for our purpose. 
First of all the two sequences $\{x_k\}$ and $\{y_k\}$
converge by (4).
Let
$$
x=\lim_{k\to\infty}x_k,\ \ y=\lim_{k\to\infty}y_k\
\mbox{ and }\ \xi=H(x),\ \ \eta=H(y).$$

By (5), $x$ and $y$ belong to the closed set $C$, and hence
$\xi$ and $\eta$ to $K=H(C)$.
By (3), we have $R_\alpha^i(x)=y$, and hence $f^i(\xi)=\eta$.
For any $k\geq 1$ (not just for $k\geq n$), define
$$
x'_k=H_k^{-1}(\xi)=H^{(k+1)}(x).
$$
By Lemma \ref{l46}, there is $m$ such that $(f_{m-1}^i)'(\xi)=(f^i)'(\xi)$.
Notice that $f_{m-1}^i$ can also be written as 
$$f_{m-1}^i=H_mR_{\alpha_m}^i H_m^{-1},$$
since $R_{\alpha_m}^i$ commutes with $h_m$.

Define 
$$\mbox{$y'_m=R_{\alpha_m}^ix'_m$ and  
$y'_k=h_{k+1}\cdots h_m(y'_m)$ for $k\leq m$.}
$$
Then we have
$$
(f_{m-1}^i)'(\xi)=(h_1)'(y'_1)\cdots(h_m)'(y_m')\,\cdot\,(h_1)'(x'_1)^{-1}
\cdots (h_m)'(x'_m)^{-1}.$$
By Lemma \ref{l43} and (4), we have 
$$
x'_k\sim_k x_k\ \mbox{ and }\ y'_k\sim_k y_k,
$$
for $k\geq n$. This shows that
$$h_n'(y_n')/h_n'(x_n')^{-1}=\lambda$$
by (1),
and for $k>n$  
\begin{equation}\label{e52}
h_k'(x'_k)=h_k'(y_k').
\end{equation}
by (3).
On the other hand we have for $k<n$,
$$
x'_k\sim_k x_n\sim_{k}y_n\sim_k y'_k.$$
This shows (\ref{e52}) for $k<n$.
The proof  that $(f^i)'(\xi)=\lambda$ is now complete.

\medskip
Now we shall construct $x_k$ and $y_k$ for $k\geq n$.

\smallskip
\noindent
{\sc Case 1 $k=n$}:
Consider a point of density of $C$ for the measure $\mu$.
Then for any $\epsilon>0$,
one can find $n\in \N$ and an interval $J$ bounded by two
consecutive fixed points of $h_n$ (a fundamental domain of $R_{Q_n}$)
such that 
$$
\mu(J\cap C)>(1-3^{-1}\epsilon)\mu(J).
$$
Then we have
$$
J\cap X_n=[x_n]_n\cup [y_n]_n,$$
where $x_n$ (resp.\ $y_n$) is the midpoint of
$[x_n]_n$ (resp.\ $[y_n]_n$), and $[x_n]_n\subset I_n^-$
and $[y_n]_n\subset I_n^+$. If $n$ is big enough, then
$\mu(J\cap X_n)$ is nearly equal to $\mu(J)$ and one
may assume
$$
\mu([x_n]_{n}\cap C>(1-2^{-1}\epsilon)\mu([x_n]_n)
\ \mbox{ and }\
\mu([y_n]_{n}\cap C>(1-2^{-1}\epsilon)\mu([y_n]_n).
$$

\smallskip
\noindent
{\sc Case 2 $k=n+1$}:
Let us call the interval $[jq_{n+1}^{-1},(j+1)q_{n+1}^{-1}]$ for some
$1\leq j\leq q_{n+1}$ a $q_{n+1}$-interval.
Now $[x_n]_n$ and $[y_n]_n$ are partitioned into $q_{n+1}$-intervals
and one can find a $q_{n+1}$-interval $J_1$ (resp.\ $J_2$)
contained in $[x_n]_n$ (resp.\ $[y_n]_n$) such that
$$\mu(J_\nu\cap C)>(1-\epsilon)\mu(J_\nu),\ \ \nu=1,2.$$

Then there is $1\leq i\leq q_{n+1}$ such that $R_{\alpha_{n+1}}^i(J_1)=J_2$.
Now since $$\mu(J_\nu\cap X_{n+1})>(1-\delta'_{n+1})\mu(J_\nu),
$$
and since $$(1-\epsilon)(1-\delta'_{n+1})>2/3,
$$
where $\delta'_{n+1}$ is the constant given by (\ref{e}),
either there are more than 2/3 portion of components $[z]_{n+1}$ among
all the components  of $J_\nu\cap X_{n+1}\cap I_{n+1}^+$ which satisfies
\begin{equation}\label{e53}
m([z]_{n+1}\cap C)\geq(1-\epsilon)(1-\delta'_{n+1})m([z]_{n+1}),
\end{equation}
 or else among all the components of 
$J_\nu\cap X_{n+1}\cap I_{n+1}^-$, for each $\nu=1,2$.
That is, we can find a component $[x_{n+1}]_{n+1}$ (resp.\
$[y_{n+1}]_{n+1}$) in $J_1$ (resp.\ $J_2$) such that 
$$R_{\alpha_{n+1}}^i(x_{n+1})
=y_{n+1}$$
which satisfies (\ref{e53}). Here we have chosen 
$x_{n+1}$ (resp.\ $y_{n+1}$) to be the midpoint of
$[x_{n+1}]_{n+1}$ (resp.\ $[y_{n+1}]_{n+1}$).

\medskip
\noindent
{\sc Case 3} {\em higher $k$}:
Assume we get $x_k$ and $y_k$ which satisfy (2), (3), (4) and

\medskip
(5')\ \ \  $\mu([x_k]_k\cap C)>(1-\epsilon)\prod_{j=1}^k(1-\delta_j')^2\mu([x_k]_k)$

and\ \ \ $\mu([y_k]_k\cap C)>(1-\epsilon)\prod_{j=1}^k(1-\delta_j')^2\mu([x_k]_k)$.

\medskip
Since $R^i_{\alpha_k}([x_k]_k)=[y_k]_k$ and since by (\ref{e33})
$$
i\abs{\alpha_{k+1}-\alpha_{k}}<2q_{n+1}\abs{\alpha_k-\alpha}
<2q_{k}\abs{\alpha_k-\alpha}<2\delta_k^2 Q_k^{-2},$$
$R_{\alpha_{k+1}}^i$ maps vast majority of $[x_k]_k$ into $[y_k]_k$.

More precisely the conditional probability
 of the union of the components of $X_{k+1}$
completely contained in $[x_k]_k\cap R_{\alpha_{k+1}}^{-i}([y_k]_k)$,
conditioned to $X_{k+1}$,
is bigger than $(1-\delta_{k+1}')^2$, and the same is true for
$R_{\alpha_{k+1}}^i([x_k]_k)\cap [y_k]_k$).
Since
$$
(1-\epsilon)\prod_{j=1}^{k+1}(1-\delta_j')^2>2/3,$$
just as before, there are $[x_{k+1}]_{k+1}$ in $[x_k]_k$ and
$[y_{k+1}]_{k+1}$ in $[y_k]_k$ such that
$
R_{\alpha_{k+1}}^i(x_{k+1})=y_{k+1}
$ and
$$
\mu([x_{k+1}]_{k+1}\cap C)>(1-\epsilon)\prod_{j=1}^{k+1}(1-\delta'_j)^2\mu([x_{k+1}]),$$
$$
\mu([y_{k+1}]_{k+1}\cap C)>(1-\epsilon)\prod_{j=1}^{k+1}(1-\delta'_j)^2\mu([y_{k+1}]).$$

This finishes the construction of $\{x_k\}$ and $\{y_k\}$.

\section{Type ${\rm III}_\infty$}

We construct $\hat h_n$ by almost the same manner as in Section 3.
But for $n$ odd we use the slope $\lambda_1^{\pm 1}$ and for $n$ even
$\lambda_2^{\pm 1}$, where $1<\lambda_1<\lambda_2$.
They are chosen so that  $\log\lambda_1$ and
$\log\lambda_2$ are independent over $\Q$.
We just repeat the argument of Section 5, to show
$\lambda_1,\ \lambda_2\in r(f)$. 

All that need extra care is the validity of 
the argument which shows that a variant of Proposition \ref{p51}
is sufficient.
For this, first of all, we need a variant
of Lemma \ref{l46}.
It holds true, as is remarked in the proof there.
We also need the fact that
$(f^i)'$ is locally constant on $\Xi\cap f^{-i}(\Xi)$.
To establish this, notice that $\log(f^i_n)'$ converges uniformly
to $\log(f^i)'$ and that $(f^i_n)'$ is locally constant on
$\Xi\cap f^{-i}(\Xi)$ if $\abs{i}\leq n$ (Compare Lemmata \ref{l45}
and \ref{l46}).
Thus for any $\xi\in\Xi\cap f^{-i}(\Xi)$,
there is a neighbourhood $U$ of $\xi$ and $n_0\in\N$
such that if $\eta\in U\cap \Xi
\cap f^{-i}(\Xi)$,
\begin{enumerate}
\item
$\abs{\log(f_{n+1}^i)'(\eta)-\log(f_n^i)'(\eta)}<\log\lambda_1$
if $n\geq n_0$ and
\item
$(f_{n_0}^i)'(\eta)=(f_{n_0}^i)'(\xi)$.
\end{enumerate}

Then we have $(f_{n}^i)'(\eta)=(f_{n_0}^i)'(\eta)$ for any $n\geq n_0$.
The same is true for $\xi$.
This shows $(f^i)'(\eta)=(f^i)'(\xi)$.

\section{Type ${\rm III}_0$}

We construct $\hat h_n$ starting at affine functions of slope $3^{-1}$
and $3^{3^n}$, and define the intervals $\hat I_n^-$ and $\hat I_n^+$
as in Section 3. Although $m(\hat I_n^+)\to0$ rapidly as $n\to\infty$,
we have $m(\hat h_n(\hat I_n^+))\to 2/3$ and $m(\hat h_n(\hat I_n^-))\to 1/3$. 

First of all let us show that $r(f)\subset\{0,1\}$.
Define $X_n$, $Y_n$ and $\Xi$ as in Section 4. Let $L$ be a component of
$H(X_n)$. We shall show  that if $n$ is sufficiently large
and $\xi,f^i(\xi)\in
L\cap\Xi$,
then either $(f^i)'(\xi)=1$ or $\abs{\log_3(f^i)'(\xi)}>3^{n-1}$.
This is sufficient for our purpose since $m(L\cap\Xi)>0$.

To show this, notice that there is $m$ bigger than $\abs{i}$ and $n$ such that
$f_m^i(\xi)\in H_m(Y_m)$ and $(f_m^i)'(\xi)=(f^i)'(\xi)$.
Define $x_j=H_j^{-1}(\xi)$ and $y_j=H_j^{-1}(f^i_m(\xi))$ for $j\leq m$.
Then we have
$$
\log_3(f_m^i)'(\xi)=\sum_{j=1}^{m}\log_3h_j'(y_j)-
\sum_{j=1}^{m}\log_3h_j'(x_j)$$
Notice that
$$\abs{\log_3h_j'(y_j)-\log_3h_j'(x_j)}= 0\ \mbox{ or }\ 3^j-1,
\ \mbox{ and }\
$$
$$\log_3h_j'(y_j)=\log_3h_j'(x_j)\ \mbox{ if }\ j\leq n,$$
since $\xi,f_m^i(\xi)\in L\cap\Xi$.

Let $k\in[n,m]$ be the largest integer, if any, such that
$$\log_3h_k'(y_k)\neq\log_3h_k'(x_k).$$
Then the value
$$
\abs{\log_3h_k'(y_k)-\log_3h_k'(x_k)}$$
is vastly bigger than 
$$
\sum_{j=1}^{k-1}\abs{\log_3h_j'(y_j)-\log_3h_j'(x_j)}.$$
In fact, computation shows that if $n$ is sufficiently large,
$$\abs{\log_3(f^i)'(\xi)}=\abs{\log_3(f_m^i)'(\xi)}>3^{n-1},
$$
as is required.

What is left is to show that $0\in r(f)$, since we always have
$1\in r(f)$. To show this, we follow the argument of Section 5
closely and show the following proposition.

\begin{proposition}
For any closed subset $K\subset\Xi$ of positive measure and any $n_0$, there
exist a point $\xi\in K$ and a number $i\in\Z$ such that 
$f^i(\xi)\in K$ and $(f^i)'(\xi)=3^{3^{n}-1}$ for some $n>n_0$.
\end{proposition}

\section{Type ${\rm II}_\infty$}

We begin by constructing $\hat h_n$ starting at affine maps
of slope $2^{\pm n}$. Define $\hat I_n^+$ just as before.
Notice that $m(\hat I_n^+)\to0$ rapidly as $n\to\infty$,
but that $m(\hat h_n(I_n^+))\to 1$ rapidly.
Define $I_n^+$ as the lift of $\hat I_n^+$ by the $Q_n$-fold covering map.
Denote
$$\ X^+=\bigcap_{i=1}^\infty I_i^+
\mbox{ and }\ \Xi^+=H(X^+)
.$$

Then one can show as in Lemma \ref{l42} that $m(\Xi^+)>0$.
On the other hand it is easy to show that $m(X^+)=0$.
This implies that the unique $f$-invariant measure $H_*m$ is singular to $m$.
Therefore $f$ cannot be of type ${\rm II}_1$.

On the other hand, we can show the following proposition easily.

\begin{proposition}
Whenever $\xi, f^i(\xi)\in\Xi^+$ for some $i\in\Z$,
then $(f^i)'(\xi)=1.$
\end{proposition}

This completes the proof that $r(f)=\{1\}$. Since $f$ is not
of type ${\rm II}_1$, it must be of type ${\rm II}_\infty$.


\begin{thebibliography}{99}

\bibitem[AK]{AK} D. Anosov and A. Katok, {\em New examples in smooth
ergodic theory, Ergodic diffeomorphisms,} Trans.\ Moscow Math.\ Soc.\ 
{\em 23} (1970), 1-35.

\bibitem[FS]{FS} B. Fayad and M. Saprykina, {\em Weak mixing disc and
annulus diffeomorphisms with arbitrary Liouvillean rotation number on the
boundary,} Ann.\ Sci.\ Ecole Norn.\ Sup.\ {\bf 38}(2005), no.\ 3,
339-364.

\bibitem[H]{H} M. R. Herman, {\em Sur la conjugaison diff\'erentiable
des diff\'eomorphismes du cercle a des rotations,} Publ.\ I. H. E. S. 
{\bf 49}(1979), 5-233.

\bibitem[KH]{KH}
A. Katok and B. Hasselblatt, {\em Introduction to the modern theory
of dynamical systems,} Encyclopedia of Mathematics and its Applications,
Vol.\ 54, Cambridge University Press 1995.

\bibitem[Ka]{Ka} Y. Katznelson, {\em Sigma-finite invariant measures for
	    smooth mappings of the circle,} J. d'Analyse Math.\ {\bf
	    31}
(1977), 1-18.

\bibitem[Kr]{Kr} W. Krieger, {\em On the Araki-Woods asumptotic ratio set
and non-singular transformations of a measure space,} In:
``Contributions to Ergodic Theory and Probability,'' Springer Lecture
	    Notes
{\bf 160}(1970), 158-177.

\bibitem[M1]{M1} S. Matsumoto, {\em A generic dimensional properties of the
invariant measures for circle diffeomorphisms,} Preprints.
arXiv:1101.4463.

\bibitem[M2]{M2} S. Matsumoto, {\em Dense properties of the space of
	    circle diffeomorphisms with a Liouville rotation number,}
Nonlinearity {\bf 25}(2012, 1495-1511.


\bibitem[S]{S} V. Sadovskaya, {\em Dimensional characteristics of
invariant measures for circle diffeomorphisms,}
Erg.\ Th.\ Dyn.\ Sys.\ {\bf 29}(2009), no.6, 1979-1992.

\bibitem[Y1]{Y1} J.-C. Yoccoz, {\em Conjugaison diff\'erentiable des 
diff\'eomorphismes du cercle dont le nombre de rotation v\'erifie une
conditon diophantine,}
Ann.\ Sci.\ Ecole Norm.\ Sup.\ {\bf 17}(1984), no.\ 3, 333-359.

\bibitem[Y2]{Y2} J.-C. Yoccoz, {\em Centralisateurs et conjugaison 
diff\'erentiable des diff\'eomorphismes du cercle,}
Ast\'erisque {\bf 231}(1995), 89-242.

\end{thebibliography}
\end{document}